\documentclass[12pt]{extarticle}

\usepackage[english]{babel}
\usepackage{amssymb,amsthm,amsmath,mathtools}
\usepackage{subcaption}
\usepackage{tikz}
\usepackage{graphicx}
\usetikzlibrary{calc}
\usetikzlibrary{decorations.pathreplacing}
\tikzstyle{snode}=[circle ,draw=black,fill=red,thick, inner sep=0pt ,minimum size=1.4mm]
\tikzstyle{bnode}=[circle ,draw=black,fill=black,thick, inner sep=0pt ,minimum size=1.4mm]
\usepackage{textcomp}
\usepackage{ifthen}
\usepackage[normalem]{ulem}

\usepackage{enumitem}
\setlist[enumerate]{left=0pt,label=(\roman*),
ref=(\roman*),font=\normalfont,topsep=-1ex,parsep=-.3ex,partopsep=0pt}

\usepackage{hyperref}
\hypersetup{
    colorlinks=true,
    pdftitle={On the independence number in subcubic graphs},
    pdfauthor={J. Harant, I. Schiermeyer},
    bookmarks=true,
    pdfpagemode=UseOutlines,
}
\usepackage[capitalise]{cleveref}

\usepackage[left=2.5cm,right=2.5cm,top=3cm,bottom=3cm]{geometry} 
\newtheorem{theorem}{Theorem}
\newtheorem{coro}{Corollary}
\newtheorem{lemma}{Lemma}

\newtheorem{proposition}{Proposition}

\title{{\bf On the independence number in subcubic graphs}}

\author{Jochen Harant$^1$ and Ingo Schiermeyer$^{2,3}$  \\~\\
\small $^1$  Technical University of  Ilmenau, 98693 Ilmenau, 
\small Germany\\
\small $^2$ AGH University of Krakow,
\small 
al. Mickiewicza 30, 30-059 Krak\'ow, 
\small Poland\\
\small $^3$ TU Bergakademie Freiberg,
\small 09596 Freiberg,
\small Germany\\
}

\begin{document}

\maketitle 

\begin{abstract}
For a connected subcubic graph 
$G\neq K_1$ let $V_i(G) = \{v \in V(G) ~|~ d_G(v)=i\}$ for $1 \leq i \leq 3.$ Given $c_1, c_2, c_ 3 \in \mathbb{R}^+$ and $ d \in \mathbb{R}$, 
 we show several results of type $\alpha(G) \geq c_1|V_1(G)| + c_2|V_2(G)| + c_3|V_3(G)| - d.$ We  also derive  classes of graphs $G$ showing sharpness of these  lower bounds on the independence number $\alpha(G)$ of $G$. 
\\[2mm]

\end{abstract}

\section{Introduction and Results}\label{I}

We consider simple, finite, and undirected graphs $G$ only, where  
 $V(G)$ and $E(G)$ denote the vertex set and the edge set of $G$. If $v\in V(G)$, then  $N_G(v)$ and $d_G(v)=|N_G(v)|$ are the neighbourhood and the degree of $v$ in $G$.  A graph $G$ with maximum degree 
$\Delta(G)$ at most $3$ is called a {\it subcubic} graph.  For terminology and notation not defined here 
we refer to \cite{West}.

A subset $I$ of $V(G)$ is called {\it independent} if the subgraph of $G$ induced by $I$ is edgeless. 
The {\it independence number $\alpha(G)$} of $G$ is the largest cardinality among all 
independent sets of $G.$ An independent set $I$ with cardinality $|I|= \alpha(G)$ is called a
{\it maximum independent set.} Finding a maximum independent set in a given graph $G$ is known as the 
{\it Maximum Independent Set Problem}. This problem is NP-hard and it remains 
NP-hard for subcubic graphs \cite{HLLM}. This tells us the following:

"The Maximum Independent Set Problem for general graphs is not harder than the Maximum Independent Set Problem for subcubic graphs."

This message gives the motivation to study this problem in subcubic graphs. For triangle-free subcubic graphs and other subclasses 
of subcubic graphs several results have been shown for the Maximum Independent Set Problem \cite{BLS,BRRS,FL,HT,HLLM,HLR,LMR,S}.

For a  subcubic graph $G$ without isolated vertices let $V_i(G) = \{v \in V(G)~ | ~d_G(v)=i\}$ for $1 \leq i \leq 3.$ Given $c_1, c_2, c_ 3 \in \mathbb{R}^+$ and $ d \in \mathbb{R},$ 
 we  show several results of type 

\begin{equation}\label{eq1}
\alpha(G) \geq c_1|V_1(G)| + c_2|V_2(G)| + c_3|V_3(G)| - d
\end{equation}

for connected subcubic graphs $G$. 

Clearly, (\ref{eq1}) holds  for all subcubic graphs $G$ without isolated vertices if and only if (\ref{eq1}) holds  for all connected subcubic graphs $G\ne K_1$. 

We restrict our consideration to the class ${\cal G}$ of connected subcubic graphs $G\neq K_4$ on at least three vertices.

We  also derive subclasses of  ${\cal G}$ for which (\ref{eq1}) is sharp. 

An immediate consequence of the famous and well known Theorem of Brooks (\cite{B}) on the relationship between the maximum degree of a graph and its chromatic number is the following Corollary \ref{Brooks}. 
\begin{coro}\label{Brooks}
If $G\in {\cal G}$, then
\begin{equation}\label{eq2}
\alpha(G) \geq \frac{1}{3}|V_1(G)| + \frac{1}{3}|V_2(G)| + \frac{1}{3}|V_3(G)|.\end{equation}
\end{coro}

A milestone in the theory on independence in graphs is the proof of the Caro - Wei - Inequality  (\cite{C,Wei})
$$ \alpha(G)\ge \sum\limits_{v\in V(G)}\frac{1}{d_G(v)+1}$$
for an arbitrary graph $G$. 

If $G$ is a graph such that $\Delta(G)\le 4$, then Brause, Randerath, Rautenbach, and Schiermeyer \cite{BRRS} improved  the result of Caro and Wei to

$$ \alpha(G)\ge \sum\limits_{v\in V(G)}\frac{2}{d_G(v)+\omega(v)+1},$$
where $\omega(v)$ denotes the maximum order of a clique of $G$ that contains the vertex $v$. 

If $G\in {\cal G}$, then $\omega(v)=2$ if $v\in V_1(G)$ and $\omega(v)\le 3$ otherwise. Thus, Corollary \ref{BRRS} follows. 
\begin{coro}\label{BRRS}~\\
If $G\in {\cal G}$, 
then 
\begin{equation}\label{eq3}
\alpha(G)\ge \frac{1}{2}|V_1(G)|+\frac{1}{3}|V_2(G)|+\frac{2}{7}|V_3(G)|.
\end{equation} 
\end{coro}

Recently, Kelly and Postle \cite{KP} proved  the following remarkable Theorem \ref{KP}.

\begin{theorem}\label{KP}~\\
If $G$ is a graph and $g:V(G) \rightarrow R$ such that $g(v)\le \frac{2}{2d_G(v)+1}$ for each $v\in V(G)$ and $\sum\limits_{v\in K}g(v)\le 1$ for each clique $K\subseteq V(G)$, then $\alpha(G)\ge \sum\limits_{v\in V(G)}g(v)$. 
\end{theorem} 

 Corollary \ref{TKP} is a simple consequence of Theorem \ref{KP}, however, its proof in Section \ref{proof} shows that, for  $G\in {\cal G}$ on at least $5$ vertices, the inequalities (\ref{eq4}) and (\ref{eq5}) are the best possible  results of the form $\alpha(G)\ge c_1|V_1(G)|+c_2|V_2(G)|+c_3|V_3(G)|$ which follow from Theorem \ref{KP}.

\begin{coro}\label{TKP}~\\
Let $G\in {\cal G}$  on at least $5$ vertices. Then 
\begin{equation}\label{eq4}
 \alpha(G)\ge \frac{2}{3}|V_1(G)|+\frac{5}{14}|V_2(G)|+\frac{2}{7}|V_3(G)|
\end{equation} and 
\begin{equation}\label{eq5}
\alpha(G)\ge \frac{3}{5}|V_1(G)|+\frac{2}{5}|V_2(G)|+\frac{1}{5}|V_3(G)|.\end{equation}
\end{coro}

Note that (\ref{eq4}) is stronger than (\ref{eq3}).  
As consequences of our forthcoming main results Theorem \ref{Main} and Theorem \ref{c1}, we  present new lower bounds on $\alpha(G)$ for $G\in {\cal G}$ and also improve the inequalities (\ref{eq2}), (\ref{eq4}), and (\ref{eq5}).  

Assume $c_1,c_2,c_3\in \mathbb{R}^{+}$, $d\in \mathbb{R}$, and that $G\in {\cal G}$. If we wish to show a result of type   $\alpha(G)\ge  c_1|V_1(G)|+c_2|V_2(G)|+c_3|V_3(G)|-d$, then we should keep the following fact in mind. 
  If $e=xy \in E(G)$ with $x,y\in V_2(G)\cup V_3(G)$, then $G'=G-e$ is subcubic (possibly no longer connected)  and $\alpha(G')\ge \alpha(G)$. Thus, the expression $c_1|V_1(G)|+c_2|V_2(G)|+c_3|V_3(G)|$ should not decrease if  $G$ is replaced with $G'$. This is the case if $ c_1\ge c_2\ge c_3$.   
  
  Proposition \ref{coeff} contains some more information about $c_1$, $c_2$, and $c_3$. Its proof  is given in Chapter \ref{proof}. 
  
 \begin{proposition}\label{coeff}~\\
If an inequality $\alpha(G)\ge c_1|V_1(G)|+c_2|V_2(G)|+c_3|V_3(G)|-d$ is fulfilled  for all $G\in {\cal G}$, then \\
(i) $c_3\le \frac{1}{3}$, \\
(ii) if $c_3= \frac{1}{3}$, then $c_2\le \frac{1}{3}$, \\
(iii) if $c_2= c_3= \frac{1}{3}$, then $c_1\le \frac{2}{3}$, and\\
(iv) $c_2\le \frac{1}{2}$. 
\end{proposition} 

We also want to find an infinite subclass  of ${\cal G}$ such that the desired inequality  \\$\alpha(G)\ge  c_1|V_1(G)|+c_2|V_2(G)|+c_3|V_3(G)|-d$ becomes tight if $G$ is chosen from this class. Therefore, two infinite classes ${\cal A}$ and ${\cal B}$ of graphs are considered. It turns out that \\ $c_3=c_3(c_1,c_2)$ and $d=d(c_1,c_2)$ as linear functions of $c_1$ and $c_2$ can be eliminated if \\$\alpha(G)\ge  c_1|V_1(G)|+c_2|V_2(G)|+c_3|V_3(G)|-d$ is tight for all graphs in ${\cal A}$ or in ${\cal B}$. Then we show that the remaining pair $(c_1,c_2)$ can be  chosen such that \\$\alpha(G)\ge  c_1|V_1(G)|+c_2|V_2(G)|+c_3(c_1,c_2)|V_3(G)|-d(c_1,c_2)$ is fulfilled for all  $G\ {\cal G}$.

Let ${\cal A}\subset {\cal G}$  be defined as follows. \\
(i) $ K_3 \in {\cal A}$.\\
(ii) If $G'\in {\cal A}$ and $x \in V_2(G')$, then the graph $G$ obtained from the disjoint union of $G'$ and $K_3$ by adding an edge between $x$ and a vertex $u\in V(K_3)$ belongs to ${\cal A}$.

Clearly, $V_2(K_3)\neq \emptyset$ and $V_2(G')\neq \emptyset$ implies $V_2(G)\neq \emptyset$ in step (ii) of the definition of ${\cal A}$, hence, ${\cal A}$ is infinite. 

Moreover, if $G\in {\cal A}$, then $V(G)=V_2(G)\cup V_3(G)$ and there is an integer $k\ge 1$ such that $|V(G)|=3k$,  $|V_2(G)|=k+2$, $|V_3(G)|=2k-2$, and, by  the mentioned Theorem of Brooks, it follows $\alpha(G)= k$.

If an inequality $\alpha(G)\ge  c_1|V_1(G)|+c_2|V_2(G)|+c_3|V_3(G)|-d$ is tight for all graphs $G$ in ${\cal A}$, then
$k=c_2(k+2)+c_3(2k-2)-d$ for all $k\ge 1$. 
Thus, $k+1=c_2(k+3)+c_3(2k)-d$ if $k$ is replaced with $k+1$. It follow $1=c_2+2c_3$ and $d=3c_2-1$. Hence, the inequality   
\begin{equation}\label{A}
\alpha(G)\ge  c_1|V_1(G)|+c_2|V_2(G)|+\frac{1-c_2}{2}|V_3(G)|-(3c_2-1) ~is ~tight ~for ~all ~G\in {\cal A} .
\end{equation}
Let ${\cal B}\subset {\cal G}$ be the infinite set of graphs  obtained from a path $P_n$ on $n\ge 2$ vertices by adding a pending edge to every vertex of $P_n$. 

 Clearly, if
 $G\in {\cal B}$, then $|V_1(G)|=|V_3(G)|+2$,  $|V_2(G)|=2$, and  $V_1(G)$ is a maximum independent set of $G$, i.e. $\alpha(G)=|V_1(G)|$. 
 
Similarly as above, it follows (\ref{B}). 
\begin{equation}\label{B}
\alpha(G)\ge  c_1|V_1(G)|+c_2|V_2(G)|+(1-c_1)|V_3(G)|-(2c_1+2c_2-2) ~is ~tight ~for ~all ~G\in {\cal B} .
\end{equation} 
Let $S_{\cal A}=\{(c_1,c_2)\in \mathbb{R}^2\}$ such that  the inequality (\ref{A}) is true for all $G\in {\cal G}$ and \\
$S_{\cal B}=\{(c_1,c_2)\in \mathbb{R}^2\}$ such that (\ref{B}) is true for all $G\in {\cal G}$.

Eventually, if an inequality $\alpha(G)\ge  c_1|V_1(G)|+c_2|V_2(G)|+c_3|V_3(G)|-d$ is tight for all graphs $G$ in ${\cal A}\cup {\cal B}$, then $c_3=\frac{1-c_2}{2}=1-c_1$ and $d=3c_2-1=2c_1+2c_2-2$ implying $c_2=2c_1-1$. 

Thus, if $S_{{\cal A}\cup {\cal B}}=\{(c_1,2c_1-1)\in S_{\cal A}\cap S_{\cal B}\}$ and if $(c_1,2c_1-1)\in S_{{\cal A}\cup {\cal B}}$, then
\begin{equation}\label{AB}
\alpha(G)\ge  c_1|V_1(G)|+(2c_1-1)|V_2(G)|+(1-c_1)|V_3(G)|-(6c_1-4)
\end{equation}
 is tight for all $G\in {\cal A}\cup {\cal B}$ and is true for all $G\in {\cal G}$.

Theorem \ref{Main} shows that $(\frac{2}{3},\frac{1}{3})\in S_{{\cal A}\cup {\cal B}}$.

\begin{theorem}\label{Main}~\\
 If  $G\in {\cal G}$,  
then 
\begin{equation}\label{eq9}
\alpha(G)\ge  \frac{2}{3}|V_1(G)|+\frac{1}{3}|V_2(G)|+\frac{1}{3}|V_3(G)|
\end{equation}
and equality holds in {\em (\ref{eq9})} if 
 $G\in {\cal A}\cup {\cal B}$.\\
 \end{theorem}

By Proposition \ref{coeff}, (\ref{eq9}) is the best possible improvement of (\ref{Brooks}).

We remark that $\alpha(G)\ge \frac{7}{12}|V_1(G)|+\frac{5}{12}|V_2(G)|+\frac{1}{3}|V_3(G)|$ for all triangle-free graphs 
 $G\in {\cal G}$ follows from a result of Griggs  \cite{griggs}  and that this inequality is tight for exactly three triangle-free graphs 
 $G\in {\cal G}$.

For the proof of the forthcoming Theorem \ref{c1} we need a similar result as presented in Lemma \ref{Griggs II}. The proof of Lemma \ref{Griggs II} is given in Chapter \ref{proof}.

Let ${\cal N}$ contain the four graphs $P_4$, disjoint copies of $C_5$ and  $K_2$ connected by an edge, $C_7$, and two disjoint copies of $C_5$ connected by an edge.

 \begin{lemma}\label{Griggs II}~\\
    If   $G$ is a connected, triangle-free, and subcubic graph on at least $2$ vertices,  
then 
\begin{equation}\label{eq10}
\alpha(G)\ge  \frac{4}{7}|V_1(G)|+\frac{3}{7}|V_2(G)|+\frac{2}{7}|V_3(G)|
\end{equation}
unless $G\in \{K_2,C_5\}$ when 
\begin{equation}\label{eq11}
\alpha(G)\ge  \frac{4}{7}|V_1(G)|+\frac{3}{7}|V_2(G)|+\frac{2}{7}|V_3(G)|-\frac{1}{7}.
\end{equation}
Moreover, equality holds in {\em (\ref{eq10})} if and only if
 $G\in {\cal N}$.\\
 \end{lemma}
 
 Note that the assertion of Lemma \ref{Griggs II} is much weaker than that ones of the forthcoming Theorem \ref{c1} and  Corollary \ref{5/7}  because $\frac{5}{7}>\frac{4}{7}$ and in Lemma \ref{Griggs II} it is assumed that $G$ is triangle free. However, Lemma \ref{Griggs II} is needed for the proof of Theorem \ref{c1}.
 
 Theorem \ref{c1} shows that $(c_1,\frac{3}{7})\in S_{\cal B}$ if $\frac{5}{7}\le c_1\le 1$. 
 
 \begin{theorem}\label{c1}~\\
    If $\frac{5}{7}\le c_1 \le 1$  and  $G\in {\cal G}$, then
    \begin{equation}\label{eq12}
   \alpha(G)\ge  c_1|V_1(G)|+\frac{3}{7}|V_2(G)|+(1-c_1)|V_3(G)|-(2c_1-\frac{8}{7}),
   \end{equation}
  where {\em (\ref{eq12})} is tight for $G\in {\cal B}$.\\
  Moreover, if $c_1=\frac{5}{7}$, then {\em (\ref{eq12})} is also tight for $G\in {\cal A}$.
  
  \end{theorem}

If $c=\frac{5}{7}$, then by Theorem \ref{c1} it follows Corollary \ref{5/7}.

\begin{coro}\label{5/7}
\begin{equation}\label{eq13} 
\alpha(G)\ge  \frac{5}{7}|V_1(G)|+\frac{3}{7}|V_2(G)|+\frac{2}{7}|V_3(G)|-\frac{2}{7}
\end{equation}
for all $G\in {\cal G}$, where equality holds if $G\in {\cal A}\cup {\cal B}$.
\end{coro}  

Note that (\ref{eq13}) is stronger than (\ref{eq4}) and also stronger than (\ref{eq5}).

By (\ref{A}), (\ref{B}), and (\ref{AB}), it follows that $S_{\cal A}$, $S_{\cal B}$, and $S_{{\cal A}\cup {\cal B}}$ are convex sets.

Thus, Theorem \ref{Main} and Theorem \ref{c1} imply  Corollary \ref{end}.

\begin{coro}\label{end}
$(c_1,2c_1-1) \in S_{{\cal A}\cup {\cal B}}$ for $\frac{1}{3}\le c_1 \le \frac{5}{7}$ and $(c_1,\frac{3}{7}) \in S_{\cal B}$ for $\frac{5}{7}\le c_1 \le 1$.
\end{coro}

We conclude with an open question, whether  Corollary \ref{end} possibly remains true if $\frac{5}{7}$ is replaced with $\frac{3}{4}$ and $\frac{3}{7}$ is replaced with $\frac{1}{2}$. Note that in this case the value $c_2=2c_1-1=\frac{1}{2}$ would be best possible. To show that the answer to this question is yes, it is sufficient to prove that
$\alpha(G)\ge  \frac{3}{4}|V_1(G)|+\frac{1}{2}|V_2(G)|+\frac{1}{4}|V_3(G)|-\frac{1}{2}$ and
 $\alpha(G)\ge  |V_1(G)|+\frac{1}{2}|V_2(G)|-1$ are true for all $G\in {\cal G}$ because $S_{\cal B}$, and $S_{{\cal A}\cup {\cal B}}$ are convex sets.

\section{Proofs}\label{proof}

\begin{proof}[{\bf Proof of Corollary \ref{TKP}}]~\\
Corollary \ref{TKP} will follow from Theorem \ref{KP}. Therefore,
let  $c_1\ge c_2\ge c_3\ge 0$ and $c_i\le \frac{2}{2i+1}$ for $i=1,2,3$. Moreover, let   $g(v)=c_i$ if $v\in V_i(G)$ for $i=1,2,3$, hence, \\$\sum\limits_{v\in V(G)}g(v)=c_1|V_1(G)|+c_2|V_2(G)|+c_3|V_3(G)|$. 
 
If $K$ is a clique of $G$, then, since $G\in {\cal G}$ on at least $5$ vertices, it follows $|K|=2$ or $|K|=3$ and, using $c_1\ge c_2\ge c_3$, we have 
$\sum\limits_{v\in K}g(v)\le c_1+c_2$ if $|K|=2$ and $\sum\limits_{v\in K}g(v)\le 2c_2+c_3$ if $|K|=3$.

If $(c_1,c_2,c_3)\in \{(\frac{2}{3},\frac{5}{14},\frac{2}{7}),(\frac{3}{5},\frac{2}{5},\frac{1}{5}) \}$, then
 $c_1\ge c_2\ge c_3\ge 0$, $c_1\le \frac23$, $c_2\le \frac25$, $c_3\le \frac27$, $c_1+c_2\le 1$, and $2c_2+c_3\le 1 $, hence, the inequalities (\ref{eq4}) and  (\ref{eq5}) follow by Theorem \ref{KP}. \\
The set $S=\{(c_1,c_2,c_3)\in \mathbb{R}^3~|~c_1\le \frac23,~ c_2\le \frac25,~ c_3\le \frac27,~ c_1+c_2\le 1~, 2c_2+c_3\le 1\} $ is a convex set. If $P=(c_1,c_2,c_3)\in \mathbb{R}^3$ is an extreme point of $S$, then, in $P$, $3$ of the $5$  inequalities $c_1\le \frac23,~ c_2\le \frac25,~ c_3\le \frac27,~ c_1+c_2\le 1$, and $2c_2+c_3\le 1$ are tight  (note that there are ${5 \choose 3}=10$ choices for  this situation) and the remaining $2$ inequalities have to be fulfilled. It  can be checked easily, that in $4$ of these $10$ cases an extreme point can be found and they are  listed in ${\cal E}=\{ (\frac{2}{3},\frac{5}{14},\frac{2}{7}),(\frac{3}{5},\frac{2}{5},\frac{1}{5}),(\frac{2}{3},\frac{1}{3},\frac{2}{7}),(\frac{9}{14},\frac{5}{14},\frac{2}{7}) \}$. 
Moreover, if $(c_1,c_2,c_3) \in {\cal E}$, then $c_1\ge c_2\ge c_3\ge 0$.
Note that $(\frac{2}{3},\frac{1}{3},\frac{2}{7})$ and $(\frac{9}{14},\frac{5}{14},\frac{2}{7})$ lead to weaker results than inequality (\ref{eq4}). 
\\Thus, if the values $g(v)$ for $v\in V(G)=V_1(G)\cup V_2(G)\cup V_3(G)$  are specifically selected as $g(v)=c_i$ if $v\in V_i(G)$ for $i=1,2,3$, then, in this case, Corollary \ref{TKP} presents the two best possible  results of the form  $\alpha(G)\ge c_1|V_1(G)|+c_2|V_2(G)|+c_3|V_3(G)|$ that can be obtained as consequences of Theorem \ref{KP} for $G\in {\cal G}$ on at least $5$ vertices.
\end{proof}

\begin{proof}[{\bf Proof of Proposition \ref{coeff}.}]~\\
To see (i), let $\varepsilon >0$ and assume that $\alpha(G)\ge c_1|V_1(G)|+c_2|V_2(G)|+c_3|V_3(G)|-d$ is fulfilled  for all  $G\in {\cal G}$, where $c_3=\frac{1}{3}+\varepsilon$.\\
Let $G'$ be a connected and $3$-regular graph on $n> |\frac{d}{3\varepsilon}|$ vertices and $G$ be obtained from $G$ by {\em blowing up} every vertex $v\in V(G')$ to a triangle such that the former edges of $G'$ form a matching in $G$. This operation is sometimes also referred to as {\em truncating} the  vertices of $G'$. Then $G$ is also connected, $3$-regular, has $3n$ vertices, and contains $n$ mutually vertex disjoint triangles, hence, $\alpha(G)\le n=\frac{1}{3}|V_3(G)|$. Clearly, $G\neq K_4$ and  because of the previously  mentioned Theorem of Brooks the chromatic number of $G$ is $3$, thus, $\alpha(G)= n=\frac{1}{3}|V_3(G)|$.\\
It follows $\frac{1}{3}|V_3(G)|\ge (\frac{1}{3}+\varepsilon)|V_3(G)|-d$, hence, $3n=|V_3(G)|\le \frac{d}{\varepsilon}$, a contradiction.\\
The proofs of (ii), (iii), and (iv) are very similar to that one of (i) and follow by the same principle. \\   
To see (ii), 
consider $G'\in {\cal A}$ on $3k$ vertices, where $k$ is large enough. Using $V_1(G')=\emptyset$, $|V_2(G')|=k+2$, $|V_3(G')|=2k-2$, and $\alpha(G')=k$, it follows 
 $\alpha(G')=\frac{1}{3}|V_2(G')|+\frac{1}{3}|V_3(G')|$.\\
For the proof of (iii), let again $G'\in {\cal A}$ on $3k$ vertices ($k$ large enough) and
 $G$ be obtained from $G'$ by adding to each vertex $x\in V_2(G')$ a disjoint copy of $K_2$ and connecting $x$ with a vertex of this $K_2$. Then  $|V_1(G)|=k+2$, $\alpha(G)=\alpha(G')+|V_1(G)|=2k+2$, $|V_2(G)|=|V_2(G')|=k+2$, and $|V_3(G)|=|V_2(G')|+|V_3(G')|=3k$, hence,
$\alpha(G)=\frac{2}{3}|V_1(G)|+\frac{1}{3}|V_2(G)|+\frac{1}{3}|V_3(G)|$.\\
Eventually, to see (iv), let $G$ be a long cycle. \end{proof}

The proofs of Theorem \ref{Main}, Lemma \ref{Griggs II}, and Theorem \ref{c1} are by induction on the number of vertices of the graph $G$.
  
Let $K_{1,3}+e$ and $K_4-e$ denote the graphs obtained from $K_{1,3}$ by adding an edge $e$ and  obtained from $K_4$ by removing an edge $e$, respectively.

\begin{proof}[{\bf Proof of Theorem \ref{Main}}]~\\
If $G\in {\cal A}$ and $|V(G)|=3k$, then recall that  $V_0(G)=V_1(G)=\emptyset$, $V(G)=V_2(G)\cup V_3(G)$, $|V_2(G)|=k+2$, $|V_3(G)|=2(k-1)$, and  $\alpha(G)= k$.  Thus,   equality holds in (\ref{eq9}) for $G\in {\cal A}$. 

If
 $G\in {\cal B}$, then $|V_1(G)|=|V_3(G)|+2$,  $|V_2(G)|=2$, and   $\alpha(G)=|V_1(G)|$ imply that equality holds in (\ref{eq9}) also for $G\in {\cal B}$.

Let  $f_G(v)=\frac{2}{3}$, $f_G(v)=\frac{1}{3}$, and $f_G(v)=\frac{1}{3}$ if  $v\in V_1(G)$, $v\in V_2(G)$, and $v\in V_3(G)$, respectively. We   use 
$\frac{2}{3}|V_1(G)|+\frac{1}{3}|V_2(G)|+\frac{1}{3}|V_3(G)|=\sum\limits_{v\in V(G)}f_G(v)$.
Note that $\alpha(G)\ge \sum\limits_{v\in V(G)}f_G(v)$ is equivalent to (\ref{eq9}).\\  
 
If $V_1(G)=\emptyset$, then Theorem \ref{Main} follows from Corollary \ref{Brooks}.

Thus, we may assume that $V_1(G)\neq \emptyset$.

For the induction base we show that Theorem \ref{Main} is true if $3\le |V(G)|\le 4$.\\
If $|V(G)|=3$, then $G=P_3$ or $G=K_3$ and Theorem \ref{Main} is true in these cases. \\If $|V(G)|=4$, then $G\in \{P_4,K_{1,3},K_{1,3}+e, K_4-e\}$. Also in these cases Theorem \ref{Main} holds.

In the sequel we assume  $|V(G)|\ge 5$.

Let $u\in V_1(G) $ and $N_G(u)=\{w\}$.  It follows $w\in V_2(G)\cup V_3(G)$.

{\em Case 1: $w\in V_2(G)$.}

Let $G'=G-\{u,w\}$ and $N_G(w)=\{u,x\}$. \\
It follows $|V(G')|\ge  3$.\\

Clearly,  $G'$ is connected, $x\in V_2(G)\cup V_3(G)$, and   $\alpha(G)=\alpha(G')+1$.\\
 It follows $f_{G'}(x)-f_G(x)\ge 0$ and, by induction, \\
$\alpha(G)= \alpha(G')+1\ge \sum\limits_{v\in V(G')}f_{G'}(v)+1=\sum\limits_{v\in V(G)}f_{G}(v) -f_G(u)-f_G(w)+f_{G'}(x)-f_G(x) +1$\\$
\ge \sum\limits_{v\in V(G)}f_{G}(v) -\frac{2}{3}-\frac{1}{3} +1= \sum\limits_{v\in V(G)}f_{G}(v)$.\\

{\em Case 2}: $w\in V_3(G)$.

Let $N_G(w)=\{u,x_1,x_2\}$ and $G'$ be obtained from $G$ by removing $\{u,w\}$ and, if $x_1x_2\notin E(G)$, by adding the edge $x_1x_2$. It follows $|V(G')|\ge 2$ and that $G'$ is connected. 

If $|V(G')|= 2$, then $G=K_{1,3}$ or $G=K_{1,3}+e$ and it follows \\
$\alpha(G)=3  \ge \frac{7}{3}=\sum\limits_{v\in V(G)}f_{G}(v) $ or
$\alpha(G)=2 \ge \frac{5}{3}=\sum\limits_{v\in V(G)}f_{G}(v) $, respectively.

We may assume  $|V(G')|\ge 3$.

Note that a maximum independent set of $G'$ is  an independent set of $G-\{u,w\}$, hence, $\alpha(G)\ge  \alpha(G')+1$.

It follows $f_{G'}(x_1)-f_{G}(x_1)\ge 0$ and $f_{G'}(x_2)-f_{G}(x_2)\ge 0$ and, by induction, \\
$\alpha(G)\ge  \alpha(G')+1\ge \sum\limits_{v\in V(G')}f_{G'}(v) +1$\\
$\ge \sum\limits_{v\in V(G)}f_{G}(v)  -\frac{2}{3}-\frac{1}{3}+1= \sum\limits_{v\in V(G)}f_{G}(v) $.
\end{proof}

\begin{proof}[{\bf Proof of Lemma \ref{Griggs II}}]~\\
We use that $G'$ is triangle-free for every subgraph $G'$ of $G$.

The statement of Theorem \ref{Griggs II} is true if $G\in \{K_2,C_5\}$. Hence, we assume $|V(G)|\ge 3$ and $G\neq C_5$. Theorem \ref{Griggs II} is proved if (\ref{eq10}) holds in this case.

Let  now $f_G(v)=\frac{4}{7}$, $f_G(v)=\frac{3}{7}$, and $f_G(v)=\frac{2}{7}$ if  $v\in V_1(G)$, $v\in V_2(G)$, and $v\in V_3(G)$, respectively. 
Note that $\alpha(G)\ge \sum\limits_{v\in V(G)}f_G(v)$ is equivalent to (\ref{eq10}). 

For the induction base we show that (\ref{eq10}) is fulfilled if  $3\le |V(G)|\le 4$.\\
If $|V(G)|=3$, then $G=P_3$ and $\alpha(G)> \sum\limits_{v\in V(G)}f_G(v)$.\\
If $|V(G)|=4$, then $G=P_4$ or $G=K_{1,3}$.
If $G=P_4$, then  $\alpha(G)= \sum\limits_{v\in V(G)}f_G(v)$ and $P_4\in {\cal N}$ while $\alpha(G)> \sum\limits_{v\in V(G)}f_G(v)$ for $G=K_{1,3}$.

In the sequel we assume $|V(G)|\ge 5$ and $G\neq C_5$.

{\em Case 1:} $V_1(G)\neq \emptyset$.\\
Let $u\in V_1(G)$ and $N_G(u)=\{w\}$.

{\em Case 1.1:} $u \in V_2(G)$.

Let $G'=G-\{u,w\}$ and $N_G(w)=\{u,x\}$.\\
Clearly, $G'$ is connected, $|V(G')|\ge 3$, and $\alpha(G)=\alpha(G')+1$.\\
It follows   $f_{G'}(x)-f_G(x)=\frac{1}{7}$, and, by induction,\\
$\alpha(G)=\alpha(G')+1\ge  \sum\limits_{v\in V(G')}f_{G'}(v)-\frac{1}{7}+1$\\
$= \sum\limits_{v\in V(G)}f_{G}(v) -\frac{1}{7}+1- \frac{4}{7}-\frac{3}{7}+ \frac{1}{7}= \sum\limits_{v\in V(G)}f_{G}(v)$. If $\alpha(G')+1=  \sum\limits_{v\in V(G')}f_{G'}(v)-\frac{1}{7}$, then equality holds, hence, $G'=C_5$,  (\ref{eq10}) is tight for $G$, and $G\in {\cal N}$ in this case.

{\em Case 1.2:} $u \in V_3(G)$.

Let $G'=G-\{u,w\}$ and $N_G(w)=\{u,x_1,x_2\}$. \\
If $G'$ is connected, then $|V(G')|\ge 3$, $f_{G'}(x_1)-f_G(x_2)=f_{G'}(x_2)-f_G(x_2)=\frac{1}{7}$, and \\
$\alpha(G)=\alpha(G')+1\ge  \sum\limits_{v\in V(G')}f_{G'}(v)-\frac{1}{7}+1$\\
$= \sum\limits_{v\in V(G)}f_{G}(v) -\frac{1}{7}+1- \frac{4}{7}-\frac{2}{7}+ 2\cdot \frac{1}{7}> \sum\limits_{v\in V(G)}f_{G}(v)$.\\
If $G'$ is not connected, then $G'$ has two components $G_1$ and $G_2$ with $x_1\in V(G_1)$, $x_2\in V(G_2)$, and $x_1x_2\notin E(G)$.\\
Note that $|V(G')|\ge 3$.\\
If $|V(G_1)|=1$ and $|V(G_2)|\ge 2$, then let $G''=G'-\{x_1\}$ and\\
$\alpha(G)=\alpha(G'')+2\ge  \sum\limits_{v\in V(G'')}f_{G''}(v)-\frac{1}{7}+2$\\
$= \sum\limits_{v\in V(G)}f_{G}(v) -\frac{1}{7}+2- 2\cdot \frac{4}{7}-\frac{2}{7}+ \frac{1}{7}> \sum\limits_{v\in V(G)}f_{G}(v)$.\\
If $|V(G_1)|\ge 2$ and $|V(G_2)|\ge 2$, then \\
$\alpha(G)=\alpha(G')+1\ge  \sum\limits_{v\in V(G')}f_{G'}(v)-\frac{2}{7}+1$\\
$= \sum\limits_{v\in V(G)}f_{G}(v) -\frac{2}{7}+1- \frac{4}{7}-\frac{2}{7}+ 2\cdot \frac{1}{7}> \sum\limits_{v\in V(G)}f_{G}(v)$.\\

{\em Case 2:} $V_1(G)=V_3(G)=\emptyset$.

Then $G$ is a cycle on at least $4$ vertices  and $G\neq C_5$. It can be seen easily that equality holds in (\ref{eq10}) for $G=C_7$, hence $C_7 \in {\cal N}$, and $\alpha(G)> \sum\limits_{v\in V(G)}f_G(v)$ in all other cases.

{\em Case 3:} $V_1(G)=\emptyset$ and $V_3(G)\neq \emptyset$.

{\em Case 3.1:}  $G$ contains a bridge $e\in E(G)$ such that one component $G_1$ of $G-\{e\}$ is  a $C_5$.

Let $G_2=G-V(G_1)$. Clearly, $V(G_2)\neq \emptyset$ and $\alpha(G)=\alpha(G_2)+2$.
\\
Since $V_1(G)=\emptyset$ we have $G_2\notin \{K_2,P_3\}$, and, because $G$ is triangle-free, $G_2\neq K_3$. It follows $|V(G_2|\ge 4$.\\
If $G_2=C_5$, then  it can be seen easily that equality holds in (\ref{eq10}) for $G$ and $G\in {\cal N}$.\\
Otherwise, it follows $\alpha(G_2)\ge  \sum\limits_{v\in V(G_2)}f_{G_2}(v)$, hence,\\
$\alpha(G)=\alpha(G_2)+2\ge  \sum\limits_{v\in V(G_2)}f_{G_2}(v)+2= \sum\limits_{v\in V(G)}f_{G}(v) +2- 4\cdot \frac{3}{7}-2\cdot \frac{2}{7}+ 2\cdot \frac{1}{7}> \sum\limits_{v\in V(G)}f_{G}(v)$ and we are done in {\em Case 3.1}.

Now let $w\in V_3(G)$.

{\em Case 3.2:} Not {\em Case 3.1} and $G-\{w\}$ has a component $G_1$ such that $G_1=C_5$.

Let $G_2=G-V(G_1)$. Clearly, $V(G_2)\neq \emptyset$ because $G$ is triangle-free.  Then $w$ has two neighbours in $V(G_1)$ and one neighbour in $V(G_2)$. Hence, $G_2\neq C_5$.\\
Again because $V_1(G)=\emptyset$ it follows $G\notin \{K_2,P_3\}$, and, becaus $G$ is triangle-free, $|V(G_2)|\ge 4$. \\
Note that $\alpha(G)=\alpha(G')+2$ for $G'=G-V(G_1)$ and $|V(G')|\ge 5$.

It follows\\
$\alpha(G)=\alpha(G')+2\ge \sum\limits_{v\in V(G')}f_{G'}(v)+2=\sum\limits_{v\in V(G)}f_{G}(v)- 3\cdot \frac{3}{7}-2\cdot \frac{2}{7}+  \frac{2}{7}+2> \sum\limits_{v\in V(G)}f_{G}(v)$.

{\em Case 3.3:} Not {\em Case 3.1} and not {\em Case 3.2}.

Let  $G'=G-\{w\}$, and $G_1,...,G_p$ be   the components of $G'$ ($1\le p\le 3$). Since $V_1(G)=\emptyset$ and $G$ is triangle-free, it follows $|V(G_q)|\ge 4$ and $G_q\neq C_5$ for $q=1,...,p$.

Thus,\\
$\alpha(G)\ge \sum\limits_{q=1}^p\alpha(G_q)\ge  \sum\limits_{q=1}^p\sum\limits_{v\in V(G_q)}f_{G_q}(v)= \sum\limits_{v\in V(G)}f_{G}(v) - \frac{2}{7}+3\cdot \frac{1}{7}> \sum\limits_{v\in V(G)}f_{G}(v)$
\end{proof}

\begin{proof}[{\bf Proof of Theorem \ref{c1}}]~\\
If
 $G\in {\cal B}$, then $|V_1(G)|=|V_3(G)|+2$,  $|V_2(G)|=2$, and   $\alpha(G)=|V_1(G)|$ imply that equality holds in (\ref{eq12}) for $G\in {\cal B}$.

The proof of Theorem \ref{c1} is organized as follows. 

Without mentioning in each case, we use that $G\in {\cal G}$.\\
For the induction base we show that the assertion of Theorem \ref{eq12} holds if  $3\le |V(G)|\le 5$.\\ 
Next we assume  that $V_1(G)=\emptyset$ ({\em Case 1}).\\
Then let  $G$ be triangle-free ({\em Case 1.1}) and in this case we use Lemma \ref{Griggs II}. \\
It is impossible that three triangles of $G$ share a common vertex of $G$ because $\Delta(G)\le 3$ would imply $G=K_4$.\\
 If  two triangles of $G$ share a common vertex, then $G$ contains $K_4-e$ as a subgraph induced by the vertices of these two triangles ({\em Case 1.2}). Otherwise, all triangles of $G$ are mutually vertex disjoint ({\em Case 1.3}). \\
Eventually we consider the case that $V_1(G)\neq \emptyset$ ({\em Case 2}).

Let  $f_G(v)=c_1$, $f_G(v)=\frac{3}{7}$, and $f_G(v)=(1-c_1)$ if  $v\in V_1(G)$, $v\in V_2(G)$, and $v\in V_3(G)$, respectively. \\Recall that 
 $\frac{5}{7}\le c_1\le 1$ and that $\alpha(G)\ge \sum\limits_{v\in V(G)}f_G(v)-(2c_1-\frac{8}{7})$ is equivalent to (\ref{eq12}).
 
For the induction base we show that Theorem \ref{eq12} is true if $3\le |V(G)|\le 5$. \\If $3\le V(G)|\le 4$, then $G\in \{P_3,K_3,P_4, K_{1,3}, K_{1,3}+e,K_4-e\}$ and it can be seen easily that (\ref{eq12}) is fulfilled for these graphs.\\
Now let $|V(G)|=5$. Clearly, $|V_1(G)|\le 3$. \\If $|V_1(G)|= 3$, then $G$ is obtained from $K_{1,3}$ by subdividing an edge by a vertex and (\ref{eq12}) is also true.\\
If $|V_1(G)|= 2$, then $V_1(G)$ does not dominate $V(G)$, hence, \\$\alpha (G)\ge 3> 2c_1+3\cdot \frac{3}{7}-(2c_1-\frac{8}{7})\ge \sum\limits_{v\in V(G)}f_G(v)-(2c_1-\frac{8}{7})$.\\
If $|V_1(G)|= 1$, then,  because the number of vertices of odd degree is even, $|V_2(G)|=3$ and $|V_3(G)|=1$ or $|V_2(G)|=1$ and $|V_3(G)|=3$. It follows that $G$ is isomorphic in the first case  to a $C_4$ with an additional  pending edge ($\alpha(G)=3$) and in the second case to  a $K_4-e$ with an additional pending edge ($\alpha(G)=2$). In both cases (\ref{eq12}) holds.

In the sequel we assume $|V(G)|\ge 6$.\\It is important to remark that  $c_1-\frac{3}{7}>\frac{3}{7}-(1-c_1)=c_1-\frac{4}{7}$ is frequently used.
 
{\em Case 1: $V_1(G)=\emptyset$.}

{\em Case 1.1: $G$ is triangle-free.}

We apply Lemma  \ref{Griggs II}, obtain 
$\alpha(G)\ge \frac{3}{7}|V_2(G)|+\frac{2}{7}|V_3(G)|-\frac{1}{7}$\\
$\ge \frac{3}{7}|V_2(G)|+(1-c_1)|V_3(G)|-(2c_1-\frac{8}{7})$, and we are done in {\em Case 1.1}.

{\em Case 1.2: $G$ contains $H=K_4-e$ on vertex set $\{u_1,u_2,w_1,w_2\}$  as a subgraph such that \\$e=w_1w_2\notin E(G)$.}

Let $G'$ be obtained from $G$ by removing $\{u_1,u_2\}$ and by adding the edge $w_1w_2$. Then $G'$ is connected and $|V(G')|\ge 4$. If $I'$ is a maximum independent set of $G'$ then $|I'\cap \{w_1,w_2\}|\le 1$. If $w_1\notin I$, then $(I'\setminus \{w_2\})\cup \{u_1,u_2\}$ is an independent set of $G$.\\
By induction, \\
$\alpha(G)\ge \alpha(G')+1\ge \sum\limits_{v\in V(G')}f_{G'}(v)-(2c_1-\frac{8}{7})+1$\\$ =\sum\limits_{v\in V(G)}f_{G}(v)-(2c_1-\frac{8}{7})-f_G(u_1)-f_G(u_2)+f_{G'}(w_2)-f_G(w_1)+f_{G'}(w_2)-f_G(w_2)+1$\\$\ge \sum\limits_{v\in V(G)}f_{G}(v)-(2c_1-\frac{8}{7})-2(1-c_1)+2(c_1-\frac{4}{7})+1$\\$=\sum\limits_{v\in V(G)}f_{G}(v)-(2c_1-\frac{8}{7})+4c_1-\frac{15}{7}> \sum\limits_{v\in V(G)}f_{G}(v)-(2c_1-\frac{8}{7})$.

{\em Case 1.3:  $G$ contains triangles but all of them are mutually vertex disjoint.}

 Let  $w_1,w_2,w_3$ be the vertices of a triangle of $G$.

{\em Case 1.3.1:  $w_1,w_2\in V_2(G), w_3\in V_3(G)$.} 

 Let $N_G(w_3)=\{w_1,w_2,u_3\}$ and $G'=G-\{w_1,w_2,w_3\}$.

It follows  that $G'$ is connected and $|V(G')|\ge 3$.\\
If  $I'$ is a maximum independent set of $G'$, then $I=I'\cup \{w_1\}$ is an independent set of $G$.

Because $u_3\in V_2(G)\cup V_3(G)$,  $f_{G'}(u_3)-f_G(u_3)\ge c_1-\frac{4}{7}$ and \\
$\alpha(G)\ge  \alpha(G')+1\ge \sum\limits_{v\in V(G')}f_{G'}(v)-(2c_1-\frac{8}{7})+1$\\$= \sum\limits_{v\in V(G)}f_{G}(v)-(2c_1-\frac{8}{7})-\frac{6}{7}-(1-c_1)+c_1-\frac{4}{7}+1\ge \sum\limits_{v\in V(G)}f_{G}(v)-(2c_1-\frac{8}{7})$.

{\em Case 1.3.2:  $w_1\in V_2(G), w_2,w_3\in V_3(G)$.} 

Let  $N_G(w_2)=\{w_1,w_3,u_2\}$ and $N_G(w_3)=\{w_1,w_2,u_3\}$. Clearly, $u_2\neq u_3$ ({\em Case 1.3}). \\
Let  $G'$ be obtained from $G$ by removing $\{w_1,w_2,w_3\}$ and adding the edge $u_2u_3$ if $u_2u_3\notin E(G)$. Then again $|V(G')|\ge 3$, $\alpha(G)\ge  \alpha(G')+1$, and $G'$ is connected.

Because $u_2,u_3\in V_2(G)\cup V_3(G)$,  $f_{G'}(u_2)-f_G(u_2)\ge 0$ and $f_{G'}(u_3)-f_G(u_3)\ge 0$ and \\
$\alpha(G)\ge \alpha(G')+1\ge \sum\limits_{v\in V(G')}f_{G'}(v)-(2c_1-\frac{8}{7})+1$\\$
\ge \sum\limits_{v\in V(G)}f_{G}(v)-(2c_1-\frac{8}{7})-\frac{3}{7}-2(1-c_1)+1\ge \sum\limits_{v\in V(G)}f_{G}(v)-(2c_1-\frac{8}{7})$ .

{\em Case 1.3.3:  $w_1, w_2,w_3\in V_3(G)$.} 
   
Let  $N_G(w_1)=\{w_2,w_3,u_1\}$, $N_G(w_2)=\{w_1,w_3,u_2\}$, and $N_G(w_3)=\{w_1,w_2,u_3\}$. Clearly, $u_1,u_2,$ and $u_3$ are mutually distinct ({\em Case 1.3}). \\  
If $u_1,u_2,$ and $u_3$ are the vertices of a triangle of $G$, then $V(G)=\{w_1,w_2,w_3,u_1,u_2,u_3 \}$ and \\$\alpha(G)=2> 6(1-c_1)-(2c_1-\frac{8}{7})=\sum\limits_{v\in V(G)}f_{G}(v)-(2c_1-\frac{8}{7})$.

Thus, we may assume that $u_1u_2\notin E(G)$ and let $G'$ be obtained from $G$ by removing $\{w_1, w_2,w_3\}$ and adding the edge $u_1u_2$. If $I'$ is a maximum independent set of $G'$, then $I'\cup \{w_1\}$  is an  independent set of $G$ if $u_1\notin I'$ and $I'\cup \{w_2\}$  is an  independent set of $G$ if $u_2\notin I'$, hence, $\alpha(G)\ge  \alpha(G')+1$. \\
$G'$ has at most two components, each component has at least three vertices because \\$V_1(G)=\emptyset$,   and  $\alpha(G')\ge \sum\limits_{v\in V(G')}f_{G'}(v)-2(2c_1-\frac{8}{7})$.

Note that $f_{G'}(u_1)=f_G(u_1)$, $f_{G'}(u_2)=f_G(u_2)$, and $f_{G'}(u_3)-f_G(u_3)\ge c_1-\frac{4}{7}$. It follows \\
$\alpha(G)\ge \alpha(G')+1\ge \sum\limits_{v\in V(G')}f_{G'}(v)-2(2c_1-\frac{8}{7})+1
\ge \sum\limits_{v\in V(G)}f_{G}(v)-2(2c_1-\frac{8}{7})-3(1-c_1)+c_1-\frac{4}{7}+1$\\$\ge  \sum\limits_{v\in V(G)}f_{G}(v)-(2c_1-\frac{8}{7})$.

{\em Case 2: $V_1(G)\neq \emptyset$.}

Let $u\in V_1(G) $ and $N_G(u)=\{w\}$.  It follows $w\in V_2(G)\cup V_3(G)$ since $|V(G)|\ge 6$.

{\em Case 2.1: $w\in V_2(G)$.}

Let $G'=G-\{u,w\}$ and $N_G(w)=\{u,x\}$. \\
Note that $|V(G')|\ge 4$.\\

Clearly,  $G'$ is connected, $x\in V_2(G)\cup V_3(G)$, and   $\alpha(G)=\alpha(G')+1$.\\
 It follows $f_{G'}(x)-f_G(x)\ge c_1-\frac{4}{7}$ and \\
$\alpha(G)= \alpha(G')+1\ge \sum\limits_{v\in V(G')}f_{G'}(v)-(2c_1-\frac{8}{7})+1$\\$
\ge \sum\limits_{v\in V(G)}f_{G}(v)-(2c_1-\frac{8}{7}) -c_1-\frac{3}{7}+c_1-\frac{4}{7} +1= \sum\limits_{v\in V(G)}f_{G}(v)-(2c_1-\frac{8}{7})$.\\

{\em Case 2.2}: $w\in V_3(G)$.

Let $N_G(w)=\{u,x_1,x_2\}$ and $G'$ be obtained from $G$ by removing $\{u,w\}$ and, if $x_1x_2\notin E(G)$, by adding the edge $x_1x_2$. It follows $|V(G')|\ge 4$ and that $G'$ is connected.

Note that a maximum independent set of $G'$ is  an independent set of $G-\{u,w\}$, hence, $\alpha(G)\ge  \alpha(G')+1$.

It follows $f_{G'}(x_1)-f_{G}(x_1)\ge 0$ and $f_{G'}(x_2)-f_{G}(x_2)\ge 0$ and, by induction, \\
$\alpha(G)\ge  \alpha(G')+1\ge \sum\limits_{v\in V(G')}f_{G'}(v)-(2c_1-\frac{8}{7}) +1$\\
$\ge \sum\limits_{v\in V(G)}f_{G}(v)-(2c_1-\frac{8}{7})  -c_1-(1-c_1)+1= \sum\limits_{v\in V(G)}f_{G}(v)-(2c_1-\frac{8}{7}) $.
\end{proof}

 \section*{Declarations}
No data have been used. There are no competing interests. All authors contributed equally.

\end{document}